\documentclass[11pt]{article}

\usepackage{amsmath, amssymb, amsthm}
\usepackage{url}
\usepackage[margin=1in]{geometry}
\usepackage[hypertexnames=false]{hyperref}
\usepackage[dvipsnames]{xcolor}
\hypersetup{
colorlinks=true,
linkcolor=NavyBlue,
citecolor=ForestGreen,
urlcolor=NavyBlue,
}
\usepackage[nameinlink,noabbrev]{cleveref}

\newcommand{\E}{\mathbb{E}}
\newcommand{\Prob}{\mathbb{P}}
\newcommand{\1}{\mathbf{1}}
\newcommand{\M}{\mathsf{M}}
\newcommand{\ALG}{\mathsf{ALG}}

\title{Threshold Rules for the Classical Prophet Inequality}
\author{Jiechen Zhang}
\date{}

\begin{document}
\maketitle

\begin{abstract}
This note records a common threshold/surplus decomposition for single-threshold stopping rules in the classical prophet inequality. The same decomposition is used to certify several deterministic thresholds, including the median, half-mean, and balanced-surplus thresholds, and to give an averaged certificate for randomized thresholds distributed as the maximum.
\end{abstract}

\section{Problem setting}
Let $X_1,\ldots,X_n$ be independent nonnegative random variables, and let $\M:=\max_{1\le i\le n}X_i$ denote the prophet's payoff. The classical prophet inequality is usually attributed to Krengel and Sucheston, who proved that there exists a nonanticipating stopping rule with expected payoff at least one half of $\E[\M]$ \cite{krengel1977semiamarts,krengel1978semiamarts}. The sharp constant $1/2$ is often attributed to Garling in this early line of work. Samuel-Cahn later showed that the same guarantee can be achieved by a single deterministic threshold rule \cite{samuelcahn1984threshold}.

We use the notation
\[
F_i(x):=\Prob(X_i\le x),
\qquad
F(x):=\prod_{i=1}^n F_i(x)=\Prob(\M\le x).
\]
For a deterministic threshold $\tau\ge 0$, define
\[
p(\tau):=\Prob(\M>\tau),
\qquad
R(\tau):=\sum_{i=1}^n \E[(X_i-\tau)^+].
\]
Throughout, we first assume that ties occur with probability zero; the general case follows by independent infinitesimal tie-breaking perturbations.

\section{A common threshold decomposition}
Let $\tau\ge 0$ be a threshold, possibly random, independent of $X_1,\ldots,X_n$. The algorithm stops at the first index $i$ such that $X_i>\tau$. Define
\[
I_i:=\1\{X_i>\tau\}\prod_{j<i}\1\{X_j\le \tau\},
\qquad
\ALG:=\sum_{i=1}^n X_iI_i.
\]
Since $I_i=1$ implies $X_i>\tau$, we have $X_iI_i=\tau I_i+(X_i-\tau)^+I_i$. Summing over $i$ and using $\sum_i I_i=\1\{\M>\tau\}$ gives
\begin{equation}
\E[\ALG]
=
\underbrace{\E[\tau\1\{\M>\tau\}]}_{\text{threshold part}}
+
\underbrace{\E\left[\sum_i (X_i-\tau)^+I_i\right]}_{\text{surplus part}}.
\label{eq:decomposition}
\end{equation}

This threshold/surplus decomposition reveals the central tradeoff: the threshold part benefits from a larger value of $\tau$ on the event $\{\M>\tau\}$, while the surplus part depends on how much value remains above $\tau$ and on the probability of reaching the corresponding item. All threshold rules below use the same decomposition; they differ only in how the surplus term is certified.

\section{Deterministic thresholds}
Assume in this paragraph that $\tau$ is deterministic. Conditioning on whether the algorithm reaches item $i$, independence gives
\[
\E[\ALG]
=
p(\tau)\tau+
\sum_{i=1}^n
\left(\prod_{j<i}F_j(\tau)\right)
\E[(X_i-\tau)^+].
\]
Since
\[
\prod_{j<i}F_j(\tau)
\ge
\prod_{j=1}^n F_j(\tau)
=
F(\tau)
=
1-p(\tau),
\]
we obtain the deterministic surplus certificate
\begin{equation}
\E[\ALG]
\ge
p(\tau)\tau+(1-p(\tau))R(\tau).
\label{eq:det-certificate}
\end{equation}
Moreover,
\[
R(\tau)
=
\sum_{i=1}^n\E[(X_i-\tau)^+]
\ge
\E[(\M-\tau)^+]
\ge
\E[\M]-\tau.
\]
Therefore
\[
\E[\ALG]
\ge
p(\tau)\tau+(1-p(\tau))(\E[\M]-\tau).
\]
Equivalently,
\begin{equation}
\E[\ALG]
\ge
\frac12\E[\M]
+
(1-2p(\tau))\left(\frac12\E[\M]-\tau\right).
\label{eq:median-halfmean-certificate}
\end{equation}

\subsection{Median and half-mean thresholds}
Equation \eqref{eq:median-halfmean-certificate} immediately gives two standard choices. If $\tau_{\mathrm{med}}$ is a median threshold of $\M$ \cite{samuelcahn1984threshold}, chosen so that
\[
p(\tau_{\mathrm{med}})
=
\Prob(\M>\tau_{\mathrm{med}})
=
\frac12,
\]
then the first factor in \eqref{eq:median-halfmean-certificate} vanishes. Likewise, for the half-mean threshold \cite{wittmann1995prophet,kleinbergweinberg2012matroid},
\[
\tau=\frac12\E[\M],
\]
the second factor vanishes. Either threshold satisfies $\E[\ALG]\ge \frac12\E[\M]$.

More generally, every deterministic threshold between $\tau_{\mathrm{med}}$ and $\frac12\E[\M]$ also works. Indeed, since $p(\tau)$ is nonincreasing in $\tau$, the two factors $1-2p(\tau)$ and $\frac12\E[\M]-\tau$ have the same sign throughout this interval. Therefore
\[
(1-2p(\tau))\left(\frac12\E[\M]-\tau\right)\ge 0,
\]
and hence $\E[\ALG]\ge \frac12\E[\M]$.

\subsection{Balanced surplus threshold}
Another useful deterministic threshold \cite{samuelcahn1984threshold} is the unique fixed point $\tau^\star$ of
\[
\tau^\star
=
R(\tau^\star)
=
\sum_{i=1}^n \E[(X_i-\tau^\star)^+].
\]
Indeed, $R$ is continuous and nonincreasing, with $R(0)=\sum_{i=1}^n \E[X_i]$ and $R(\tau)\to 0$ as $\tau\to\infty$. Moreover, $R(\tau)-\tau$ is strictly decreasing: if $0\le s<t$, then $R(t)-t\le R(s)-t<R(s)-s$. Hence $R(\tau)-\tau$ has exactly one zero.

Using \eqref{eq:det-certificate},
\[
\E[\ALG]
\ge
p(\tau^\star)\tau^\star+(1-p(\tau^\star))R(\tau^\star)
=
\tau^\star.
\]
On the other hand,
\[
\E[\M]
\le
\tau^\star+\E[(\M-\tau^\star)^+]
\le
\tau^\star+\sum_{i=1}^n\E[(X_i-\tau^\star)^+]
=
2\tau^\star.
\]
Thus $\E[\ALG]\ge \tau^\star\ge \frac12\E[\M]$.

The same argument shows that every
\[
\tau\in\left[\frac12\E[\M],\tau^\star\right]
\]
also works. Since $R$ is nonincreasing and $R(\tau^\star)=\tau^\star$, every such $\tau$ satisfies $R(\tau)\ge R(\tau^\star)=\tau^\star\ge \tau$. Therefore \eqref{eq:det-certificate} gives
\[
\E[\ALG]
\ge
p(\tau)\tau+(1-p(\tau))R(\tau)
\ge
p(\tau)\tau+(1-p(\tau))\tau
=
\tau
\ge
\frac12\E[\M].
\]

\subsection{A certified interval of deterministic thresholds}
Combining the two deterministic families above, the present certificates guarantee the $1/2$ bound for every
\[
\tau\in
\left[
\min\left\{\tau_{\mathrm{med}},\frac12\E[\M]\right\},
\max\left\{\tau_{\mathrm{med}},\tau^\star\right\}
\right].
\]
That is, every threshold in this interval satisfies $\E[\ALG]\ge \frac12\E[\M]$. 

\section{Random thresholds}

\subsection{A surplus formula for random thresholds}
Let $\tau$ have CDF $G$, independent of $X_1,\ldots,X_n$. Conditional on $\tau=t$, the rule becomes the deterministic threshold-$t$ rule. For item $i$ to contribute surplus, all previous values must be at most $t$, while $(X_i-t)^+$ already enforces $X_i>t$. Thus
\[
\E[(X_i-\tau)^+I_i\mid \tau=t]
=
\E[(X_i-t)^+]\prod_{j<i}F_j(t).
\]
Using the tail identity
\[
\E[(X_i-t)^+]=\int_t^\infty (1-F_i(x))\,dx
\]
and applying Tonelli's theorem, we obtain
\begin{equation}
\begin{aligned}
\E\left[\sum_{i=1}^n (X_i-\tau)^+I_i\right]
&=
\int_0^\infty
\sum_{i=1}^n
\left(\int_t^\infty (1-F_i(x))\,dx\right)
\prod_{j<i}F_j(t)\,dG(t) \\
&=
\int_0^\infty
\sum_{i=1}^n (1-F_i(x))
\left(
\int_0^x \prod_{j<i}F_j(t)\,dG(t)
\right)dx.
\end{aligned}
\label{eq:random-surplus}
\end{equation}
This form is useful when the threshold itself is randomized.

\subsection{A randomized analogue of the median threshold}
The deterministic median threshold chooses a fixed value $\tau_{\mathrm{med}}$ such that $\Prob(\M>\tau_{\mathrm{med}})=1/2$. We now consider a randomized analogue of this idea. Let $\tau$ be an independent draw from the distribution of $\M$; equivalently, $\tau\overset{d}{=}\M$, with $\tau$ independent of the online instance. In the notation above, this means $G=F$. Under the no-ties assumption, $\Prob(\M>\tau)=1/2$, so the random threshold has the same median-type exceedance probability as the deterministic median threshold. We prove the averaged surplus bound
\begin{equation}
\E\left[\sum_{i=1}^n (X_i-\tau)^+I_i\right]
\ge
\frac12\E[(\M-\tau)^+].
\label{eq:random-surplus-bound}
\end{equation}

Using \eqref{eq:random-surplus} with $G=F$, fix $x\ge 0$. For $t\le x$,
\[
F(t)
=
\left(\prod_{j<i}F_j(t)\right)
\left(\prod_{\ell\ge i}F_\ell(t)\right)
\le
\left(\prod_{j<i}F_j(t)\right)
\left(\prod_{\ell\ge i}F_\ell(x)\right).
\]
Hence
\[
\prod_{j<i}F_j(t)
\ge
\frac{F(t)}{\prod_{\ell\ge i}F_\ell(x)}.
\]
Assuming $F$ is continuous,\footnote{The continuity assumption is used here only to evaluate the Lebesgue--Stieltjes integral. If $F$ has jumps, then, writing $\Delta F(a):=F(a)-F(a-)$ and using $F(0-)=0$ since $\M\ge0$, one has $\int_{[0,x]} F(t)\,dF(t)=\frac12 F(x)^2+\frac12\sum_{0\le a\le x}(\Delta F(a))^2\ge \frac12 F(x)^2$. Thus the possible jumps only contribute a nonnegative correction at this point.}
this gives
\[
\int_0^x \prod_{j<i}F_j(t)\,dF(t)
\ge
\frac{1}{\prod_{\ell\ge i}F_\ell(x)}
\int_0^x F(t)\,dF(t)
=
\frac{F(x)^2}{2\prod_{\ell\ge i}F_\ell(x)}.
\]
Substituting this lower bound into \eqref{eq:random-surplus} yields
\[
\E\left[\sum_{i=1}^n (X_i-\tau)^+I_i\right]
\ge
\frac12
\int_0^\infty
F(x)
\sum_{i=1}^n (1-F_i(x))\prod_{j<i}F_j(x)\,dx.
\]
The sum telescopes:
\[
\sum_{i=1}^n (1-F_i(x))\prod_{j<i}F_j(x)
=
1-F(x).
\]
Therefore
\[
\E\left[\sum_{i=1}^n (X_i-\tau)^+I_i\right]
\ge
\frac12\int_0^\infty F(x)(1-F(x))\,dx.
\]
Since $\M$ and $\tau$ are independent with common CDF $F$,
\[
\E[(\M-\tau)^+]
=
\int_0^\infty \Prob(\tau<x<\M)\,dx
=
\int_0^\infty \Prob[\tau<x]\Prob[\M>x]\,dx
=
\int_0^\infty F(x)(1-F(x))\,dx.
\]
This proves \eqref{eq:random-surplus-bound}.

Combining \eqref{eq:decomposition} and \eqref{eq:random-surplus-bound},
\[
\E[\ALG]
\ge
\E[\tau\1\{\M>\tau\}]
+
\frac12\E[(\M-\tau)^+]
=
\E[\tau\1\{\M>\tau\}]
+
\frac12\E[(\M-\tau)\1\{\M>\tau\}]
=
\frac12\E[(\M+\tau)\1\{\M>\tau\}].
\]
By exchangeability of the i.i.d. pair $(\M,\tau)$ and the no-ties assumption,
\[
\E[(\M+\tau)\1\{\M>\tau\}]
=
\E[(\M+\tau)\1\{\tau>\M\}]
=
\frac12\E[\M+\tau]
=
\E[\M].
\]
Thus $\E[\ALG]\ge \frac12\E[\M]$.

A direct proof of the same guarantee via pointwise tail dominance is given in \Cref{app:tail-dominance}.

\subsection{Other randomized thresholds}

The argument above also suggests a sufficient condition for other randomized thresholds. Let $\tau$ have CDF $G$, independent of $X_1,\ldots,X_n$. Suppose first that, for every $i$ and every $x\ge 0$,
\begin{equation}
\int_0^x \prod_{j<i}F_j(t)\,dG(t)
\ge
\frac12 F(x)\prod_{j<i}F_j(x).
\label{eq:general-random-surplus-certificate}
\end{equation}
Then \eqref{eq:random-surplus} gives
\[
\E\left[\sum_{i=1}^n (X_i-\tau)^+I_i\right]
\ge
\frac12
\int_0^\infty
F(x)
\sum_{i=1}^n (1-F_i(x))\prod_{j<i}F_j(x)\,dx
=
\frac12\int_0^\infty F(x)(1-F(x))\,dx.
\]

It remains to compare the threshold part. If, in addition,
\begin{equation}
\int_0^\infty t(1-F(t))\,dG(t)
\ge
\int_0^\infty t(1-F(t))\,dF(t),
\label{eq:general-threshold-certificate}
\end{equation}
then the two terms in the decomposition are at least the corresponding two terms in the proof above. Hence the same conclusion follows:
\[
\E[\ALG]\ge \frac12\E[\M].
\]
Thus any threshold distribution $G$ satisfying \eqref{eq:general-random-surplus-certificate} and \eqref{eq:general-threshold-certificate} is also certified by the same threshold/surplus argument.

A simple source of additional randomized thresholds comes from convexity. Let $\tau_0$ be any deterministic threshold already certified above, for example any threshold in the certified deterministic interval. For $\alpha\in[0,1]$, draw $\tau$ from the distribution of $\M$ with probability $\alpha$, and set $\tau=\tau_0$ with probability $1-\alpha$, independently of the online instance. Equivalently, the law of $\tau$ is the convex combination
\[
\alpha F+(1-\alpha)\delta_{\tau_0}.
\]
Since the expected payoff is affine in the law of the independent threshold,
\[
\E[\ALG]
=
\alpha\,\E[\ALG\mid \tau\overset{d}{=}\M]
+
(1-\alpha)\,\E[\ALG\mid \tau=\tau_0]
\ge
\frac12\E[\M].
\]

\subsection{Implementing the randomized analogue}
The randomized analogue only requires drawing an independent threshold $\tau\overset{d}{=}\M$. There are several equivalent ways to implement this.

First, if we have sample access to each distribution $D_i$, draw independent samples $S_i$ from $D_i$ for $i=1,\ldots,n$, and set
\[
\tau:=\max_i S_i.
\]
Then
\[
\tau\overset{d}{=}\max_i X_i=\M,
\]
so this is exactly the sample-max threshold rule \cite{rubinsteinwangweinberg2020samples}.

Second, if the CDFs $F_i$ are known, then the CDF of $\M$ is
\[
F(x)=\prod_{i=1}^n F_i(x).
\]
Draw $U$ uniformly from $[0,1]$ and set
\[
\tau:=F^{-1}(U),
\qquad
F^{-1}(u):=\inf\{x:F(x)\ge u\}.
\]
Then $\tau\overset{d}{=}\M$.

Finally, if an independent simulator for the full instance is available, draw $(\widetilde X_1,\ldots,\widetilde X_n)$ independently and set
\[
\tau:=\max_i \widetilde X_i.
\]
If only historical data are available, one may approximate this by drawing uniformly from empirical past maxima. The exact theorem applies whenever the resulting threshold is an independent draw from the distribution of $\M$.

\section{Takeaway}
The deterministic and randomized threshold rules above are all consequences of the same decomposition
\[
\E[\ALG]
=
\underbrace{\E[\tau\1\{\M>\tau\}]}_{\text{threshold part}}
+
\underbrace{\E\left[\sum_i (X_i-\tau)^+I_i\right]}_{\text{surplus part}}.
\]

\begin{itemize}
\item For deterministic thresholds, the surplus part is
\[
\sum_{i=1}^n
\left(\prod_{j<i}F_j(\tau)\right)
\E[(X_i-\tau)^+].
\]
The factor $\prod_{j<i}F_j(\tau)$ is the probability that the rule reaches item $i$, and it is bounded below by $F(\tau)=\prod_{j=1}^n F_j(\tau)$. This gives the pointwise surplus bound $\E\left[\sum_i (X_i-\tau)^+I_i\right]\ge F(\tau)R(\tau)$.

\item For randomized thresholds, the same surplus term is averaged over the possible threshold values. After the tail-integral rewriting, it becomes
\[
\int_0^\infty
\sum_{i=1}^n (1-F_i(x))
\left(
\int_0^x \prod_{j<i}F_j(t)\,dG(t)
\right)dx.
\]
Thus the inner quantity $\int_0^x \prod_{j<i}F_j(t)\,dG(t)$ is the randomized analogue of the deterministic reach probability: it averages the probability of reaching item $i$ over all threshold values $t\le x$. The randomized surplus certificate is therefore no longer a pointwise comparison at a single threshold value. Instead, it is a lower bound on the averaged reach term $\int_0^x \prod_{j<i}F_j(t)\,dG(t)$.
For example, one sufficient bound is
\[
\int_0^x \prod_{j<i}F_j(t)\,dG(t)
\ge
\frac12 F(x)\prod_{j<i}F_j(x).
\]
\end{itemize}

The distinction is similar in spirit to the difference between pointwise dominance and integrated dominance conditions: deterministic thresholds compare reach probabilities at one threshold value, while randomized thresholds compare their averages over the threshold distribution.

This pointwise certificate is what certifies the deterministic threshold rules above, including the median threshold $\tau_{\mathrm{med}}$, the half-mean threshold $\frac12\E[\M]$, the balanced surplus threshold $\tau^\star$, and every threshold in the deterministic interval.

The randomized analogue $\tau\overset{d}{=}\M$ is the cleanest averaged certificate, since taking $G=F$ makes the surplus term telescope and lets the threshold part close by exchangeability. The same viewpoint also certifies any threshold distribution $G$ satisfying the two sufficient conditions above, as well as convex mixtures of the law of $\M$ with any certified deterministic threshold.

Therefore, the median rule, half-mean rule, balanced surplus rule, sample-max rule, and these mixture rules are different ways of certifying the same threshold/surplus tradeoff.

\bibliographystyle{alpha}
\bibliography{ref}

\appendix

\section{Alternative proof: pointwise tail dominance}
\label{app:tail-dominance}

The averaged surplus certificate \eqref{eq:random-surplus-bound} for the randomized threshold $\tau\overset{d}{=}\M$ can also be strengthened to the following pointwise tail certificate. Under the same no-ties assumption, or, for distributions with atoms, under the same independent infinitesimal tie-breaking convention used above, whenever $\tau$ is independent of $X_1,\ldots,X_n$, for every $z\ge 0$,
\begin{equation}
\Prob(\ALG>z)\ge \frac12\Prob(\M>z).
\label{eq:tail-dominance}
\end{equation}
Integrating \eqref{eq:tail-dominance} over $z$ recovers $\E[\ALG]\ge \frac12\E[\M]$. Thus the dominance can be seen pointwise in the tail level $z$, rather than only after averaging the surplus over the threshold distribution.

In the tie-breaking version, all comparisons such as $\M>\tau$ are interpreted in the tie-broken order, while the event $\ALG>z$ refers to the original, unperturbed payoff.

Let $F$ be the CDF of $\M$, and fix $z\ge 0$. On the event $\{\tau>z\}$, if $\M>\tau$, then some online value exceeds the threshold; the algorithm accepts a value larger than $\tau$, and hence larger than $z$. Since $\M$ and $\tau$ are independent and identically distributed, exchangeability, together with the no-ties/tie-breaking convention, gives
\[
\Prob(\ALG>z,\ \tau>z)
\ge
\Prob(\M>\tau,\ \tau>z)
=
\frac12(1-F(z))^2.
\]
It remains to treat $\{\tau\le z\}$. On $\{\M>z\}$, let $I$ be the first index with $X_I>z$. To represent the independent threshold, write $\tau=\max_k S_k$, where the $S_k$ are independent draws from the corresponding distributions, independent of the online $X_k$. Define
\[
Y_X:=\max_{j<I}X_j,
\qquad
Y_S:=\max_{j<I}S_j.
\]
With the empty maximum equal to $0$, conditional on $\{I=i\}$ and $\{\tau\le z\}$, the two collections $(X_j)_{j<i}$ and $(S_j)_{j<i}$ are independent copies, with each corresponding coordinate conditioned to lie at most $z$. Hence $\Prob(Y_S\ge Y_X\mid I=i,\tau\le z)\ge \frac12$. Since $\tau=\max_k S_k$, we also have $\tau\ge Y_S$. Therefore, conditional on $\{I=i\}$ and $\{\tau\le z\}$, with probability at least $1/2$, we have $\tau\ge Y_X$. In that case the algorithm does not stop before $i$, while $X_i>z\ge\tau$, so it accepts a value larger than $z$. Therefore
\[
\Prob(\ALG>z,\ \tau\le z)
\ge
\frac12\Prob(\M>z,\ \tau\le z)
=
\frac12F(z)(1-F(z)).
\]
Adding the two cases gives
\[
\Prob(\ALG>z)
\ge
\frac12(1-F(z))^2+\frac12F(z)(1-F(z))
=
\frac12(1-F(z)),
\]
which is exactly \eqref{eq:tail-dominance}.

\medskip
\noindent\rule{\linewidth}{0.4pt}

{\small\itshape
This proof was ``found'' through a prompt to ChatGPT 5.5 Plus, which produced the idea after about four minutes of thinking. The author thanks Pranav Nuti, who later pointed out that a closely related pointwise stochastic-dominance argument appears in Section~5, page~12, of Tomer Ezra's \emph{Prophet Inequality from Samples: Is the More the Merrier?}~\hyperref[appref:ezra2024samples]{[Ezr24]}.
\par}

{\renewcommand{\refname}{References for the appendix}
}

\end{document}